\documentclass[12pt,a4paper,leqno]{amsart}
\usepackage[T1]{fontenc}
\usepackage[utf8]{inputenc}
\usepackage{lmodern}
\usepackage[a4paper,left=26mm,right=26mm,top=27mm,bottom=28mm,headheight=15pt,headsep=12pt]{geometry}
\usepackage{amsmath,amssymb,mathtools}
\usepackage{microtype}
\usepackage{cite}
\usepackage[hidelinks,unicode]{hyperref}
\hypersetup{pdftitle={Curvature estimate for shrinking gradient Ricci soliton},pdfauthor={Fengjiang Li, Jianyu Ou, Yuanyuan Qu, and Guoqiang Wu},pdfsubject={Bounded scalar curvature and an entropy-gap curvature estimate in dimension five}}
\newtheorem{theorem}{Theorem}[section]
\newtheorem{proposition}[theorem]{Proposition}
\newtheorem{lemma}[theorem]{Lemma}

\theoremstyle{remark}
\newtheorem{remark}[theorem]{Remark}
\numberwithin{equation}{section}
\DeclareMathOperator{\Ric}{Ric}
\DeclareMathOperator{\Rm}{Rm}
\DeclareMathOperator{\Vol}{Vol}
\DeclareMathOperator{\AVR}{AVR}
\DeclareMathOperator{\supp}{supp}
\DeclareMathOperator{\id}{id}
\newcommand{\R}{\mathbb R}
\newcommand{\Acal}{\mathcal A}
\newcommand{\Wcal}{\mathcal W}
\newcommand{\dd}{\,d}
\newcommand{\inner}[2]{\langle #1,#2\rangle}
\newcommand{\norm}[1]{\lVert #1\rVert}
\newcommand{\Pcal}{\mathcal P}
\newcommand{\gcan}{g_{\mathrm{can}}}
\newcommand{\gE}{g_{\mathrm E}}
\newcommand{\eps}{\varepsilon}
\allowdisplaybreaks[1]
\title[Curvature estimates for shrinking gradient Ricci solitons]{Curvature estimates for shrinking gradient Ricci solitons}
\author[F. Li]{Fengjiang Li}
\address{Mathematical Science Research Center, Chongqing University of Technology, Chongqing 400054, China}
\email{fengjiangli@cqut.edu.cn}
\author[J. Ou]{Jianyu Ou}
\address{School of Mathematical Sciences, Xiamen University, Xiamen 361005, China}
\email{oujianyu@xmu.edu.cn}
\author[Y. Qu]{Yuanyuan Qu}
\address{School of Mathematical Sciences, Shanghai Key Laboratory of Pure Mathematics and Mathematical Practice, East China Normal University, Shanghai 200241, China}
\email{52285500012@stu.ecnu.edu.cn}
\author[G. Wu]{Guoqiang Wu}
\address{School of Science, Zhejiang Sci-Tech University, Hangzhou 310018, China}
\email{gqwu@zstu.edu.cn}
\date{}

\begin{document}
\begin{abstract}
Suppose $(M^5,g,f)$ is a complete gradient shrinking Ricci soliton with bounded scalar curvature. We prove that if $\nu(g)>-\log 2$, then its Riemannian curvature is bounded. The proof uses curvature blow-up, a local Ricci estimate, and the splitting given by the rescaled potential. An entropy estimate then allows us to apply the volume-ratio gap theorem for Ricci-flat four-manifolds.
\end{abstract}
\maketitle

\section{Introduction}
A gradient Ricci soliton is a complete Riemannian manifold $(M^n,g)$ with a smooth function $f$ satisfying
\[
\Ric+\nabla^2f=\lambda g
\]
for a constant $\lambda\in\R$. The soliton is shrinking, steady, or expanding according as $\lambda>0$, $\lambda=0$, or $\lambda<0$. In the shrinking case we rescale the metric so that
\begin{equation}\label{eq:soliton-intro}
\Ric+\nabla^2f=\frac12g.
\end{equation}

Gradient Ricci solitons generate self-similar solutions to the Ricci flow by diffeomorphisms and scaling. They play a basic role in the analysis of singularities of the flow. In his study of ancient $\kappa$-solutions, Perelman \cite{Perelman2002,Perelman2003} established the structure underlying the three-dimensional singularity analysis. Together with the work of Ni--Wallach \cite{NiWallach2008} and Cao--Chen--Zhu \cite{CaoChenZhu2008}, this led to the classification of three-dimensional shrinking gradient Ricci solitons. In higher dimensions, classification results usually require additional curvature assumptions; see, for example, \cite{Brendle2013,Brendle2014,CaoChen2013,Naber2010,PetersenWylie2010,PetersenWylie2009,WuWuWylie2018,Zhang2009}.

Curvature boundedness is an important preliminary issue in these classification problems. In dimension four, Munteanu--Wang \cite{MunteanuWang2015} proved that a scalar curvature upper bound implies an estimate $|\Rm|\leq CR$, and hence bounded curvature. In dimensions five and higher, the scalar curvature does not give a direct bound for the curvature tensor. It is natural to seek conditions under which bounded scalar curvature still rules out curvature concentration at infinity.

In this paper, we prove the following entropy-gap criterion in dimension five. We use Perelman's scale-invariant entropy $\nu(g)=\inf_{\tau>0}\mu(g,\tau)$, whose definition is recalled in Section~\ref{sec:prelim}.

\begin{theorem}\label{thm:main}
Let $(M^5,g,f)$ be a fixed complete gradient shrinking Ricci soliton satisfying \eqref{eq:soliton-intro}. If
\begin{equation}\label{eq:main-hyp}
\sup_M R_g<\infty,\qquad \nu(g)>-\log 2,
\end{equation}
then its Riemannian curvature is bounded:
\[
\sup_M|\Rm_g|<\infty.
\]
\end{theorem}

We prove the theorem by contradiction. A point-selection argument, the noncollapsing theorem of Li--Wang \cite{LiWang2020}, and their backward pseudolocality theorem \cite{LiWang2024} produce a complete smooth curvature blow-up limit. Bounded scalar curvature makes this limit Ricci flat. To prove that the limit splits, we use the local Ricci estimate of Chen--Wang \cite[Theorem~1.2]{ChenWang2013}, in the higher-derivative form recorded by Bamler \cite[Proposition~3.6]{Bamler2018}, to obtain the bound
\[
|\Ric_{Q_jg}|\leq C Q_j^{-1/2}
\]
on each fixed pointed ball. Using this estimate, we show that the normalized potential converges to a function with parallel unit gradient. Thus the nonflat limit has a global splitting
\[
Y^5=\R\times X^4,\qquad \Ric_{h_X}=0.
\]

Wang's continuity result for the local $\mu$-functional under pointed smooth Cheeger--Gromov convergence, together with domain monotonicity and exhaustion, provides one way to compare the entropies; see \cite{Wang2018}. We give the needed one-sided inequality directly by transferring compactly supported test functions and using the scale invariance of $\nu$. This yields $\nu(g)\leq\nu(Y)$. A Gaussian test-function argument then gives
\[
\nu(g)\leq\nu(Y)\leq\log\AVR(Y),
\]
where $\AVR$ denotes the Euclidean-normalized asymptotic volume ratio. The product splitting preserves this ratio. The entropy assumption therefore implies
\[
\AVR(X^4)=\AVR(Y^5)>\frac12.
\]
The sharp four-dimensional Ricci-flat gap theorem of Honda--Mondino \cite[Proposition~3.7]{HondaMondino2026} forces $X$ to be Euclidean, contradicting the nonflatness of the normalized blow-up. The zeroth-order local Ricci estimate is sufficient for the Hessian estimate that produces the splitting.

The paper is organized as follows. Section~\ref{sec:prelim} recalls the soliton identities, entropy, and local estimates used below. Section~\ref{sec:blowup} constructs the Ricci-flat blow-up and proves its global real-line splitting under bounded scalar curvature. Section~\ref{sec:entropy} establishes the entropy--volume estimate and invariance of the asymptotic volume ratio under Euclidean products. In Section 5, we prove Theorem~\ref{thm:main}.

\noindent\textbf{Acknowledgement:}
The present work is a continuation of the research undertaken  by the third author as a master's student. The third author would like to thank Professor Meng Zhu and Wenqi Li for their constant encouragement and helpful discussions.

\medskip

\noindent\textbf{Disclosure on AI assistance.}
The research problem and the principal mathematical ideas were
formulated by the authors. The authors used ChatGPT to assist with
checking some computations and polishing the exposition.
 All mathematical
claims, citations, and final arguments remain the responsibility of the
authors.

\section{Preliminaries}\label{sec:prelim}
Throughout the paper, manifolds are smooth, connected, and without boundary. A shrinking gradient Ricci soliton is normalized by
\begin{equation}\label{eq:soliton}
\Ric_g+\nabla_g^2f=\frac12g.
\end{equation}
After adding a constant to $f$, we also assume
\begin{equation}\label{eq:hamilton}
R_g+|\nabla f|_g^2=f.
\end{equation}
We write $\Delta_f=\Delta-\inner{\nabla f}{\nabla\cdot}$. We use the metric $g$ to take gradients, Hessians, traces, and norms unless another metric is written. The standard soliton identities are
\begin{align}
 dR&=2\Ric(\nabla f,\cdot),\label{eq:dR}\\
 \Delta f+R&=\frac n2,\label{eq:trace}\\
 \Delta_f f&=\frac n2-f,\label{eq:driftf}\\
 \Delta_f R&=R-2|\Ric|^2.\label{eq:driftR}
\end{align}
For completeness, these identities follow directly from \eqref{eq:soliton}. Taking its trace gives \eqref{eq:trace}. The contracted Bianchi identity and the divergence formula for the Hessian give
\[
0=\frac12\nabla_jR+\nabla_j\Delta f+\Ric_{j\ell}\nabla_\ell f.
\]
Since $\nabla_j\Delta f=-\nabla_jR$, this is \eqref{eq:dR}. It follows that
\begin{align*}
d\bigl(R+|\nabla f|^2-f\bigr)
 &=2\Ric(\nabla f,\cdot)+2\nabla^2f(\nabla f,\cdot)-df\\
 &=0.
\end{align*}
Thus adding a constant to $f$ indeed gives \eqref{eq:hamilton}. Subtracting this identity from \eqref{eq:trace} proves \eqref{eq:driftf}. Finally, taking the divergence of $\nabla_jR=2\Ric_{j\ell}\nabla_\ell f$ gives
\begin{align*}
\Delta R
 &=2\nabla_j\Ric_{j\ell}\nabla_\ell f+2\Ric_{j\ell}\nabla_j\nabla_\ell f\\
 &=\inner{\nabla R}{\nabla f}+R-2|\Ric|^2,
\end{align*}
which proves \eqref{eq:driftR}. In particular,We also know the fact that every complete gradient shrinking Ricci soliton has nonnegative scalar curvature \cite{Chen2009}.

We first recall the entropy conventions. For a complete manifold $(N^n,g)$ put
\[
\Acal(N,g)=\left\{u\in C_c^\infty(N):u\geq0,\ \int_Nu^2\dd V_g=1\right\}.
\]
For $u\in\Acal(N,g)$ and $\tau>0$, define
\begin{equation}\label{eq:W-u}
\begin{split}
\Wcal(g,u,\tau)
 &=\int_N\left[\tau\bigl(4|\nabla u|_g^2+R_gu^2\bigr)-u^2\log u^2\right]\dd V_g\\
 &\quad-n-\frac n2\log(4\pi\tau),
\end{split}
\end{equation}
where $u^2\log u^2$ is interpreted as zero at $u=0$. Perelman's functionals \cite{Perelman2002} are
\begin{equation}\label{eq:munu}
\mu(g,\tau)=\inf_{u\in\Acal(N,g)}\Wcal(g,u,\tau),
\qquad \nu(g)=\inf_{\tau>0}\mu(g,\tau).
\end{equation}
Constant metric rescaling satisfies
\begin{equation}\label{eq:entropy-scaling}
\mu(ag,a\tau)=\mu(g,\tau),\qquad \nu(ag)=\nu(g)\quad(a>0).
\end{equation}
To check the scaling formula, let $a>0$, put $g_a=ag$, and set $u_a=a^{-n/4}u$. Then
\[
dV_{g_a}=a^{n/2}dV_g,\qquad
R_{g_a}=a^{-1}R_g,\qquad
|\nabla u_a|_{g_a}^2=a^{-n/2-1}|\nabla u|_g^2.
\]
Consequently $u_a\in\Acal(N,g_a)$ if and only if $u\in\Acal(N,g)$. The energy term in $\Wcal(g_a,u_a,a\tau)$ is unchanged. The logarithmic term satisfies
\[
-\int_Nu_a^2\log u_a^2\dd V_{g_a}
=-\int_Nu^2\log u^2\dd V_g+\frac n2\log a.
\]
The last summand cancels the change in $-\frac n2\log(4\pi a\tau)$. Hence $\Wcal(ag,u_a,a\tau)=\Wcal(g,u,\tau)$. Taking the two infima proves \eqref{eq:entropy-scaling}.

For a probability potential $F$, with $u^2=(4\pi\tau)^{-n/2}e^{-F}$, the equivalent expression is
\begin{equation}\label{eq:W-F}
\Wcal(g,F,\tau)=\int_N\left[\tau\bigl(|\nabla F|_g^2+R_g\bigr)+F-n\right]
\frac{e^{-F}}{(4\pi\tau)^{n/2}}\dd V_g.
\end{equation}
Indeed, on the set where $u>0$,
\[
4|\nabla u|_g^2=u^2|\nabla F|_g^2,\qquad
-\log u^2=F+\frac n2\log(4\pi\tau).
\]
Substitution in \eqref{eq:W-u} gives \eqref{eq:W-F}. We will use this expression for a noncompactly supported density only after showing that its value is the limit of values for compactly supported smooth test functions. This approximation is supplied in the proof of Lemma~\ref{lem:entropy-avr}.

\begin{proposition}[{\cite{LiWang2020}}]\label{prop:entropy}
Let $(M^n,g,f)$ be a complete shrinker satisfying \eqref{eq:soliton}--\eqref{eq:hamilton}. Then
\begin{equation}\label{eq:entropy-potential}
\nu(g)=\mu(g,1)=\log\left((4\pi)^{-n/2}\int_Me^{-f}\dd V_g\right)>-\infty.
\end{equation}
If $F=f+\nu(g)$, then $(4\pi)^{-n/2}\int_Me^{-F}\dd V_g=1$ and $\Wcal(g,F,1)=\nu(g)$.
\end{proposition}

The constant added to obtain a probability potential should be kept separate from the normalization \eqref{eq:hamilton}. If
\[
I=(4\pi)^{-n/2}\int_Me^{-f}\dd V_g,
\]
then Proposition~\ref{prop:entropy} says that $I=e^{\nu(g)}$. Thus $F=f+\log I$ satisfies
\[
(4\pi)^{-n/2}\int_Me^{-F}\dd V_g=e^{-\log I}I=1.
\]
Adding this constant changes neither $\nabla f$ nor $\nabla^2f$. We continue to use $f$ with normalization \eqref{eq:hamilton} when we differentiate the potential.

\begin{lemma}\label{lem:potential}
Let $(M^n,g,f)$ be a complete noncompact shrinker satisfying \eqref{eq:soliton}--\eqref{eq:hamilton}. For a fixed reference point $o$, there is a constant $C$ such that
\begin{equation}\label{eq:potential-growth}
\frac14\bigl(d_g(o,x)-C\bigr)_+^2\leq f(x)
\leq\frac14\bigl(d_g(o,x)+C\bigr)^2.
\end{equation}
In particular, $f$ is proper and tends to infinity along every divergent sequence. Outside a compact set, the function $\rho_f=2\sqrt f$ is smooth and satisfies
\begin{equation}\label{eq:rho-identities}
|d\rho_f|_g^2=1-\frac{4R_g}{\rho_f^2},\qquad
\nabla_g^2\rho_f=\frac1{\rho_f}\bigl(g-2\Ric_g-d\rho_f\otimes d\rho_f\bigr).
\end{equation}
Moreover, $|d\rho_f|_g\leq1$ wherever $f>0$.
\end{lemma}
\begin{proof}
The growth estimate is the standard potential estimate of Cao--Zhou \cite{CaoZhou2010}; enlarging $C$ absorbs the choice of the reference point. The positive part is denoted by $(a)_+=\max\{a,0\}$. Differentiating $\rho_f=2\sqrt f$ gives $d\rho_f=f^{-1/2}df$. Equation~\eqref{eq:hamilton} gives the gradient identity in \eqref{eq:rho-identities}. For the Hessian,
\begin{align*}
\nabla^2\rho_f
 &=f^{-1/2}\nabla^2f-\tfrac12f^{-3/2}df\otimes df\\
 &=\frac1{\rho_f}\bigl(g-2\Ric-d\rho_f\otimes d\rho_f\bigr),
\end{align*}
where we used \eqref{eq:soliton}. The gradient bound follows from $R_g\geq0$. The lower growth bound implies that $f>0$ outside a compact set, so all the differentiations above are valid there. It also implies that sublevel sets of $f$ are bounded. They are closed, and hence compact by completeness. This proves the properness assertion.
\end{proof}

We next state the local results used in the curvature blow-up. For a Ricci flow $g(t)$, write
\[
\Pcal_g(x_0,t_0;r,-T,0)
=B_{g(t_0)}(x_0,r)\times[t_0-T,t_0].
\]
Here the spatial ball is formed with the metric at the terminal time $t_0$. It is the same set of points at each time in the product. Curvature norms and scalar curvature at time $t$ are computed using $g(t)$. All such neighborhoods are understood to lie in the time domain of the flow.

\begin{theorem}[Scalar-curvature noncollapsing, {\cite{LiWang2020}}]\label{thm:noncollapse}
For a fixed complete shrinker $(M^n,g,f)$, there is a constant $\kappa>0$ such that
\[
R_g\leq r^{-2}\ \text{on }B_g(x,r)
\quad\Longrightarrow\quad
\Vol_g B_g(x,r)\geq\kappa r^n.
\]
The statement is invariant under constant metric rescaling.
\end{theorem}

Only a scalar curvature bound is used in this noncollapsing statement. No bound for the full curvature tensor is required. For example, if $R_g\leq S$ on all of $M$, then every ball with $r\leq S^{-1/2}$ satisfies its hypothesis. To check its scaling invariance, put $g_a=ag$. A ball of radius $r$ for $g_a$ is the ball of radius $r/\sqrt a$ for $g$, and
\[
\Vol_{g_a}B_{g_a}(x,r)=a^{n/2}\Vol_g B_g(x,r/\sqrt a).
\]
The scalar curvature assumption transforms into $R_g\leq(r/\sqrt a)^{-2}$, and the same $\kappa$ therefore works for $g_a$. The assertion is also unchanged by a pullback under a diffeomorphism.

\begin{theorem}[Backward pseudolocality, {\cite[Theorem~6.21]{LiWang2024}}]\label{thm:pseudolocal}
For every $0<\alpha<1$, there is an $\eps=\eps(n,\alpha)>0$ with the following property. Let $g(t)$ be a Ricci flow induced by a complete shrinking gradient Ricci soliton. If
\[
\Vol_{g(t_0)}B_{g(t_0)}(x_0,r)\geq\alpha r^n,
\qquad |\Rm_{g(t_0)}|\leq(\alpha r)^{-2}
\quad\text{on }B_{g(t_0)}(x_0,r),
\]
then
\[
|\Rm_{g(t)}|\leq(\eps r)^{-2}
\quad\text{on }\Pcal_g\bigl(x_0,t_0;(1-\alpha)r,-(\eps r)^2,0\bigr).
\]
\end{theorem}

In the preceding theorem, both hypotheses concern only the terminal ball. The conclusion gives control on a smaller ball for a definite interval before $t_0$. The constants depend only on $n$ and $\alpha$, so they do not change when we apply the theorem at different points of the same rescaled flow. In Section~\ref{sec:blowup} we will apply it at every point of a fixed terminal ball. This gives the curvature bound on a whole spatial region, which is needed for the next estimate.

\begin{proposition}[Local Ricci estimate, {\cite[Proposition~3.6]{Bamler2018}}]\label{prop:local-ricci}
Let $g(t)$, $t\in[-2,0]$, be a Ricci flow on an $n$-dimensional manifold, let $(x_0,t_0)\in M\times[-1,0]$, and let $0<r<1$. Suppose $B_{g(t_0)}(x_0,r)$ is relatively compact and
\[
|\Rm|<r^{-2},\qquad |R|\leq\eta<1
\quad\text{on }\Pcal_g(x_0,t_0;r,-r^2,0),
\]
where $\eta>0$. Then, for every integer $\ell\geq0$,
\begin{equation}\label{eq:local-ricci}
|\nabla^\ell\Ric|(x_0,t_0)
\leq C_\ell(n)\eta^{1/2}r^{-\ell-1}.
\end{equation}
\end{proposition}

The zeroth-order estimate was proved earlier by Chen and Wang in \emph{On the conditions to extend Ricci flow~(III)} \cite[Theorem~1.2]{ChenWang2013}. Their theorem concerns a complete flow with $|\Rm|\leq2$ on
\[
B_{g(0)}(x_0,1)\times[-1/8,0]
\]
and gives
\[
\sup_{B_{g(0)}(x_0,1/2)\times[-1/16,0]}|\Ric|
\leq A(n)\left(
\sup_{B_{g(0)}(x_0,1)\times[-1/8,0]}|R|
\right)^{1/2}.
\]
After parabolic rescaling, this has the form
$|\Ric|(x_0,t_0)\leq C(n)r^{-1}\eta^{1/2}$.
Thus Chen--Wang's theorem gives the case $\ell=0$ needed for the splitting argument below. The higher-derivative estimates in \eqref{eq:local-ricci} are quoted from \cite[Proposition~3.6]{Bamler2018}.

The estimate \eqref{eq:local-ricci} is local and does not require the ambient manifold to be compact. In its application below, the radius $r$ is fixed in the rescaled flow, while $\eta$ tends to zero. The square-root rate is what controls the rescaled Hessian of the potential. More precisely, $\ell=0$ gives $|\Ric|\leq C_0(n)r^{-1}\eta^{1/2}$. The constant is independent of $\eta$, so this is a quantitative rate as the scalar curvature tends to zero. Smooth convergence alone would give only $|\Ric|\to0$, without a specified rate.

\section{Curvature blow-up and splitting}\label{sec:blowup}
In this section we establish the blow-up and splitting statements used in the proof of the main theorem. These statements hold in arbitrary dimension. The scalar curvature bound will be used both to obtain a Ricci-flat limit and, through a local Ricci estimate, to produce a parallel direction in that limit.

\begin{lemma}[Point selection]\label{lem:point-selection}
Let $(M,g)$ be a complete noncompact Riemannian manifold with $\sup_M|\Rm_g|=\infty$. There are points $p_j\to\infty$ and numbers $Q_j=|\Rm_g|(p_j)$ such that
\begin{equation}\label{eq:point-selection}
Q_j\geq j^4,\qquad
|\Rm_g|\leq4Q_j
\quad\text{on }B_g\left(p_j,\frac j{\sqrt{Q_j}}\right).
\end{equation}
\end{lemma}
\begin{proof}
Fix $o\in M$. Curvature is bounded on each compact subset, so we can choose $x_j$ with $d_g(o,x_j)\geq j+2$ and $|\Rm_g|(x_j)\geq j^4$. Indeed, the closed ball $\overline B_g(o,j+2)$ is compact by completeness and has bounded curvature. Unboundedness of curvature therefore allows us to choose $x_j$ outside this ball with curvature as large as required.

Start with $x_j$ as the current point. If the current point is $x$ and $Q=|\Rm_g|(x)$, stop if $|\Rm_g|\leq4Q$ on $B_g(x,j/\sqrt Q)$. Otherwise choose a point in this ball where the curvature is greater than $4Q$, and repeat.

After $k$ moves the curvature is at least $4^k|\Rm_g|(x_j)$. Consequently, counting the first move as $k=0$, the distance moved at the $k$th step is at most $2^{-k}j/\sqrt{|\Rm_g|(x_j)}$. Thus the total distance from $x_j$ is at most
\[
\frac{2j}{\sqrt{|\Rm_g|(x_j)}}\leq\frac2j.
\]
The same geometric series bounds the distance between any two late points by a tail tending to zero. An infinite sequence of moves would therefore be Cauchy. By completeness it would converge to a point of $M$, while its curvature tends to infinity. This contradicts continuity of curvature at the limit point. The process therefore terminates at a point $p_j$ satisfying \eqref{eq:point-selection}. Moreover, $d_g(o,p_j)\geq j+2-2/j$, so $p_j\to\infty$.
\end{proof}

\begin{lemma}[Parabolic control and a local Ricci estimate]\label{lem:parabolic}
Let $(M^n,g,f)$ be a fixed complete noncompact shrinker satisfying \eqref{eq:soliton} and \eqref{eq:hamilton}, with $0\leq R_g\leq S$. Let $p_j,Q_j$ satisfy \eqref{eq:point-selection}. Write
\begin{equation}\label{eq:canonical}
\begin{gathered}
\gcan(s)=(-s)\varphi_s^*g,\qquad
\frac d{ds}\varphi_s=\frac{\nabla_g f}{-s}\circ\varphi_s,
\qquad \varphi_{-1}=\id,\\
g_j(t)=Q_j\gcan(-1+t/Q_j).
\end{gathered}
\end{equation}
There are $c>0$ and $K<\infty$, independent of $j$, such that, for each fixed $L<\infty$ and all sufficiently large $j$,
\begin{equation}\label{eq:parabolic-bound}
|\Rm_{g_j(t)}|\leq K
\quad\text{on }B_{g_j(0)}(p_j,L+2)\times[-c,0].
\end{equation}
For every nonnegative integer $m$ there is also a constant $C_m$, independent of $j$ and $L$, such that, for all sufficiently large $j$ depending on $L$,
\begin{equation}\label{eq:ricci-rate}
\sup_{B_{g_j(0)}(p_j,L)}|\nabla^m\Ric_{g_j(0)}|_{g_j(0)}
\leq C_mQ_j^{-1/2}.
\end{equation}
The constants may depend on the fixed shrinker, $n$, $S$, and, in the last estimate, on $m$.
\end{lemma}
\begin{proof}
We first note that the flow in \eqref{eq:canonical} is well defined for $t<Q_j$ and satisfies $g_j(0)=Q_jg$. The potential growth estimate and $|\nabla f|^2\leq f$ show that $\nabla f$ has at most linear growth. Its time-dependent multiple $(-s)^{-1}\nabla f$ thus has a complete flow on every compact subinterval of $(-\infty,0)$. Differentiating the first formula in \eqref{eq:canonical} and using $2\nabla^2f=g-2\Ric_g$ gives
\[
\partial_s\gcan(s)=\varphi_s^*(-g+2\nabla^2f)=-2\Ric_{\gcan(s)}.
\]
Parabolic rescaling therefore makes $g_j(t)$ a Ricci flow as well.

Theorem~\ref{thm:noncollapse} gives $\kappa>0$, depending on the fixed shrinker, such that
\begin{equation}\label{eq:noncollapse-application}
\Vol_gB_g(q,r)\geq\kappa r^n
\quad\text{whenever } R_g\leq r^{-2}\text{ on }B_g(q,r).
\end{equation}
We apply Theorem~\ref{thm:pseudolocal}, the backward part of Li--Wang's pseudolocality theorem \cite[Theorem~6.21]{LiWang2024}. This theorem does not require an a priori global curvature bound.

Choose
\[
0<\alpha<\min\{\kappa,1/2\},\qquad b=(2\alpha)^{-1}.
\]
Fix $L$. If $j>L+2+b$ and $q\in B_{g_j(0)}(p_j,L+2)$, then \eqref{eq:point-selection} gives
\[
|\Rm_{g_j(0)}|\leq4=(\alpha b)^{-2}
\quad\text{on }B_{g_j(0)}(q,b).
\]
Here $B_{g_j(0)}(q,b)$ lies in $B_{g_j(0)}(p_j,j)$, so the point-selected curvature bound applies to this entire ball. Since $R_{g_j(0)}=Q_j^{-1}R_g\leq S/Q_j\leq b^{-2}$ for large $j$, the rescaled form of \eqref{eq:noncollapse-application} gives
\[
\Vol_{g_j(0)}B_{g_j(0)}(q,b)\geq\kappa b^n\geq\alpha b^n.
\]
Let $\eps=\eps(n,\alpha)>0$ be the pseudolocality constant. Backward pseudolocality gives a curvature bound $(\eps b)^{-2}$ on a backward neighborhood of $(q,0)$. In particular, this bound holds along the worldline with spatial coordinate $q$ for $-(\eps b)^2\leq t\leq0$. We apply the theorem separately at every terminal point $q\in B_{g_j(0)}(p_j,L+2)$. Thus the same bound holds on the whole fixed spatial region for this common time interval, proving \eqref{eq:parabolic-bound}. We can take $K=(\eps b)^{-2}$ and $c\leq(\eps b)^2$. The same constants work for every fixed $L$; only the lower bound on $j$ changes. We may in particular assume $c\leq1$.

The scalar curvature of a pullback metric is the pullback of the scalar curvature. Constant metric rescaling divides it by the scaling factor. With $s=-1+t/Q_j$, we have
\[
g_j(t)=(Q_j-t)\varphi_s^*g,\qquad
R_{g_j(t)}(x)=\frac{R_g(\varphi_s(x))}{Q_j-t}.
\]
The global scalar curvature bound on $g$ therefore implies, for $t\leq0$,
\begin{equation}\label{eq:scalar-rescaled}
0\leq R_{g_j(t)}\leq\frac S{Q_j-t}\leq\frac S{Q_j}.
\end{equation}
On the fixed spatial regions in \eqref{eq:parabolic-bound}, the curvature bound also bounds $\Ric$ and hence controls the variation of the metric for $-c\leq t\leq0$. More precisely, since $|\Ric|\leq C(n)K$ there, the Ricci flow equation gives
\[
e^{-2C(n)Kc}g_j(0)\leq g_j(t)\leq e^{2C(n)Kc}g_j(0).
\]
Choose $a>0$ so small that $a<1$, $a^2<c$, and $K<a^{-2}$. For $q\in B_{g_j(0)}(p_j,L)$, the terminal ball $B_{g_j(0)}(q,a)$ is contained in $B_{g_j(0)}(p_j,L+2)$ and is relatively compact by completeness. Hence its backward parabolic neighborhood of radius $a$ is contained in the region where both \eqref{eq:parabolic-bound} and \eqref{eq:scalar-rescaled} hold. The choice of $a$ depends only on $c$ and $K$, and is independent of $j$, $L$, and $q$.

We now apply Proposition~\ref{prop:local-ricci}, in the form of Bamler \cite[Proposition~3.6]{Bamler2018}, with $r=a$ and $\eta=\max\{S,1\}/Q_j$, using \eqref{eq:scalar-rescaled}. It gives
\[
|\nabla^m\Ric_{g_j(0)}|(q)
\leq A_m\left(\frac{\max\{S,1\}}{Q_j}\right)^{1/2}a^{-m-1}
\]
and proves \eqref{eq:ricci-rate}. The time-domain assumption in that proposition causes no further restriction: $g_j(t)$ is defined on all of $[-2,0]$, whereas its curvature bound is needed only on the smaller parabolic neighborhood. Also, $\max\{S,1\}/Q_j<1$ for large $j$, as required in that proposition. Absorbing the fixed factors $a^{-m-1}$ and $\max\{S,1\}^{1/2}$ into $C_m$ proves the claimed estimate uniformly at every point of the fixed ball. For $m=0$, the same estimate follows directly from Chen--Wang \cite[Theorem~1.2]{ChenWang2013} after parabolic rescaling.
\end{proof}

\begin{proposition}[A Ricci-flat blow-up limit]\label{prop:blowup}
Let $(M^n,g,f)$ be a fixed complete noncompact shrinker satisfying \eqref{eq:soliton} and \eqref{eq:hamilton}. If $0\leq R_g\leq S$ and $\sup_M|\Rm_g|=\infty$, then, for some sequence $p_j\to\infty$ and $Q_j=|\Rm_g|(p_j)\to\infty$,
\begin{equation}\label{eq:pointed-limit}
(M,Q_jg,p_j)\longrightarrow(Y^n,h,p_\infty)
\end{equation}
in the pointed smooth topology. The limit is complete and Ricci flat, and
\begin{equation}\label{eq:limit-properties}
|\Rm_h|(p_\infty)=1,\qquad \nu(Y,h)\geq\nu(M,g).
\end{equation}
\end{proposition}
\begin{proof}
Choose $p_j,Q_j$ by Lemma~\ref{lem:point-selection}, and use the rescaled flows of Lemma~\ref{lem:parabolic}. At the terminal time, \eqref{eq:point-selection} gives $|\Rm_{Q_jg}|\leq4$ on $B_{Q_jg}(p_j,j)$. For large $j$, \eqref{eq:noncollapse-application} also gives
\[
\Vol_{Q_jg}B_{Q_jg}(p_j,1)\geq\kappa,
\]
since $R_{Q_jg}\leq S/Q_j\leq1$. The curvature and volume bounds on the unit ball give, by the Cheeger--Gromov--Taylor estimate \cite{CheegerGromovTaylor1982}, a positive lower bound for the injectivity radius at $p_j$, independent of $j$. The parabolic curvature bound \eqref{eq:parabolic-bound} and local derivative estimates \cite{Shi1989} give bounds for all curvature derivatives on every fixed smaller spatial region and on a common smaller backward time interval. Pointed Hamilton compactness gives a smooth limit Ricci flow on that interval. Here we take a convergent subsequence on a fixed pointed ball, then enlarge the radius through positive integers and take a diagonal subsequence. The resulting convergence is smooth on every compact subset, including at terminal time. The limit has complete spatial slices: the controlled terminal radii tend to infinity, and the metric comparison in the proof of Lemma~\ref{lem:parabolic} holds uniformly on every fixed ball on the common backward interval.

By \eqref{eq:scalar-rescaled}, the scalar curvature of the limiting flow is identically zero. Its scalar evolution equation
\[
(\partial_t-\Delta)R=2|\Ric|^2
\]
shows that its Ricci tensor vanishes on the interior of the time interval, and hence also at terminal time by continuity. The Ricci flow equation now gives $\partial_t h(t)=0$. Thus the limit is static; write its terminal slice as $(Y,h,p_\infty)$. Smooth convergence at the basepoint gives $|\Rm_h|(p_\infty)=1$.

It remains to verify the direction of the entropy inequality. Let $u\geq0$ be smooth and compactly supported on $Y$, with $\int_Yu^2\dd V_h=1$, and let $\tau>0$. Transfer $u$ to $M$ by a pointed convergence map $\Psi_j$ whose domain contains a neighborhood of $\supp u$. Define $\widetilde u_j=u\circ\Psi_j^{-1}$ on the image of this domain and extend it by zero on $M$. The extension is smooth because $u$ vanishes near the boundary of the chosen domain. Put
\[
c_j^2=\int_M\widetilde u_j^2\dd V_{Q_jg},\qquad u_j=c_j^{-1}\widetilde u_j.
\]
Smooth convergence of the metrics on the fixed support gives $c_j\to1$. Thus $u_j\in\Acal(M,Q_jg)$ for all large $j$. The volume forms, gradient terms, and scalar curvature terms converge smoothly on this support. The entropy terms also converge, since $z^2\log z^2$ extends continuously to zero. Consequently,
\[
\Wcal(Q_jg,u_j,\tau)\longrightarrow\Wcal(h,u,\tau).
\]
Scale invariance of $\nu$ yields
\[
\nu(M,g)=\nu(M,Q_jg)\leq\Wcal(Q_jg,u_j,\tau).
\]
Passing to the limit gives $\nu(M,g)\leq\Wcal(h,u,\tau)$. Since both $u$ and $\tau$ were arbitrary, taking the infimum over this full class gives $\nu(M,g)\leq\nu(Y,h)$ and proves \eqref{eq:limit-properties}.
\end{proof}

\begin{lemma}[A global parallel direction]\label{lem:splitting}
The limit in Proposition~\ref{prop:blowup} admits a globally defined smooth function $v_\infty$ satisfying
\begin{equation}\label{eq:parallel}
|\nabla v_\infty|_h=1,\qquad \nabla_h^2v_\infty=0.
\end{equation}
Consequently,
\begin{equation}\label{eq:split}
(Y^n,h)\cong(\R\times X^{n-1},dt^2+h_X),
\end{equation}
where $X$ is complete, Ricci flat, and nonflat.
\end{lemma}
\begin{proof}
Use the normalization \eqref{eq:hamilton} and set $\rho=2\sqrt f$ outside a compact set. Lemma~\ref{lem:potential}, based on \cite{CaoZhou2010}, shows that $f$ is proper, $\rho(p_j)\to\infty$, and the identities \eqref{eq:rho-identities} hold. In particular, $\rho$ is $1$-Lipschitz. Since the selected ball has radius $j/\sqrt{Q_j}\leq1/j$ in $g$, it lies outside any fixed compact subset for all sufficiently large $j$. The Lipschitz bound therefore gives
\[
\inf_{B_g(p_j,j/\sqrt{Q_j})}\rho\geq\rho(p_j)-\frac1j\longrightarrow\infty.
\]
For each fixed $L$ and all $j>L$, the ball $B_{Q_jg}(p_j,L)$ is contained in this selected ball. Thus $\rho$ tends uniformly to infinity on every fixed pointed ball in the rescaled metrics.

Put $g_j=Q_jg$ and define
\begin{equation}\label{eq:vj}
v_j=\sqrt{Q_j}\bigl(\rho-\rho(p_j)\bigr).
\end{equation}
Then $v_j(p_j)=0$, and
\begin{equation}\label{eq:vj-gradient}
|dv_j|_{g_j}^2=|d\rho|_g^2=1-\frac{4R_g}{\rho^2}\longrightarrow1
\end{equation}
locally uniformly. To estimate its Hessian, recall that constant metric rescaling leaves the Levi--Civita connection unchanged and divides the norm of a covariant two-tensor by $Q_j$. In particular, $\Ric_{g_j}=\Ric_g$ as covariant two-tensors, while $|\Ric_g|_g=Q_j|\Ric_{g_j}|_{g_j}$. Since $|g|_g=\sqrt n$ and $|d\rho\otimes d\rho|_g=|d\rho|_g^2\leq1$, we obtain
\[
|\nabla_{g_j}^2v_j|_{g_j}
=Q_j^{-1/2}|\nabla_g^2\rho|_g
\leq\frac{\sqrt n+1+2|\Ric_g|_g}{\sqrt{Q_j}\rho}.
\]
Using the scaling of the Ricci norm, this becomes
\begin{equation}\label{eq:vj-hessian}
\begin{split}
|\nabla_{g_j}^2v_j|_{g_j}
&\leq\frac{\sqrt n+1}{\sqrt{Q_j}\rho}
+\frac{2\sqrt{Q_j}}\rho|\Ric_{g_j}|_{g_j}\\
&\leq\frac{\sqrt n+1}{\sqrt{Q_j}\rho}+\frac C\rho
\longrightarrow0,
\end{split}
\end{equation}
where the second line uses \eqref{eq:ricci-rate} with $m=0$. The convergence is uniform on every fixed pointed ball. In particular, no bound on $Q_j/f(p_j)$ is needed.

Pull back $v_j$ by the convergence maps in \eqref{eq:pointed-limit}. The estimate $|dv_j|_{g_j}\leq1$ and $v_j(p_j)=0$ bound $v_j$ on each fixed pointed ball. In a fixed coordinate chart on $Y$, smooth convergence of the metrics also bounds the first coordinate derivatives of these pulled-back functions. Their second coordinate derivatives satisfy
\[
\partial_a\partial_bv_j=(\nabla_{g_j}^2v_j)_{ab}
+(\Gamma_{g_j})^d_{ab}\partial_dv_j.
\]
The right-hand side is uniformly bounded on each compact subset. Thus the first derivatives are locally equicontinuous. Arzel\`a--Ascoli and a diagonal subsequence give local $C^1$ convergence to a function $v_\infty$ on all of $Y$, with $v_\infty(p_\infty)=0$.

Pass the Hessian equation against a compactly supported smooth function in such a chart and passing to the limit, We obtain, in the distributional sense,
\[
\partial_a\partial_bv_\infty=\Gamma^d_{ab}(h)\partial_dv_\infty.
\]
The right-hand side is continuous, so the first derivatives are $C^1$. Repeating this argument with the smooth coefficients of $h$ shows that $v_\infty$ is smooth. Hence $\nabla_h^2v_\infty=0$ classically. Local $C^1$ convergence and \eqref{eq:vj-gradient} also give $|dv_\infty|_h=1$, proving \eqref{eq:parallel}.

Let $V=\nabla v_\infty$ and $X=v_\infty^{-1}(0)$. Since $Y$ is complete and $|V|=1$, each integral curve of $V$ has speed one and cannot leave every compact subset in finite time. The flow $\Phi_t$ therefore exists for every $t\in\R$. Moreover, $\mathcal L_Vh=2\nabla_h^2v_\infty=0$, so every $\Phi_t$ is an isometry. Along this flow,
\[
\frac d{dt}v_\infty(\Phi_t(y))=|V|^2=1.
\]
Every flow line therefore meets $X$ exactly once. The map
\[
\R\times X\longrightarrow Y,\qquad(t,x)\longmapsto\Phi_t(x)
\]
is a smooth bijection, with smooth inverse
\[
y\longmapsto\bigl(v_\infty(y),\Phi_{-v_\infty(y)}(y)\bigr).
\]
Its differential sends the coordinate vector in the $\R$ direction to $V$. Tangent vectors to $X$ are perpendicular to $V$, and the flow preserves their inner products. Its pullback metric is therefore $dt^2+h_X$, proving the global isometry \eqref{eq:split}. The hypersurface $X$ is closed and totally geodesic, hence complete, and the product curvature formula gives $\Ric_{h_X}=0$. Finally, $X$ cannot be flat, because otherwise $Y$ would be flat, contrary to $|\Rm_h|(p_\infty)=1$.
\end{proof}

\section{Entropy bounds and asymptotic volume growth on Ricci-flat manifolds}\label{sec:entropy}
Throughout this section, $\omega_n$ denotes the volume of the unit ball in Euclidean $\R^n$. For a complete connected Ricci-flat $n$-manifold $(N,h)$ and $o\in N$, define
\begin{equation}\label{eq:avr}
\AVR(N,h):=\lim_{r\to\infty}\frac{\Vol_h B_h(o,r)}{\omega_n r^n}.
\end{equation}
Bishop--Gromov comparison shows that the ratio in \eqref{eq:avr} is nonincreasing and belongs to $(0,1]$. Its limit therefore exists in $[0,1]$ and is independent of $o$. Indeed, if $o'\in N$ and $D=d_h(o,o')$, then, for $r>D$,
\[
B_h(o,r-D)\subset B_h(o',r)\subset B_h(o,r+D).
\]
After division by $\omega_n r^n$, both outer volume ratios tend to the same limit. Thus changing the centre does not change the asymptotic volume ratio. The entropy $\nu(N,h)$ is understood as the infimum of the $\Wcal$-functional over positive scales and nonnegative, compactly supported smooth test functions of unit $L^2$ norm. It is allowed to equal $-\infty$.

\begin{lemma}[Entropy and asymptotic volume ratio]\label{lem:entropy-avr}
Let $(N^n,h)$ be a complete connected Ricci-flat manifold. Then
\begin{equation}\label{eq:entropy-avr}
\nu(N,h)\leq\log\AVR(N,h),
\end{equation}
where $\log0=-\infty$. In particular, if $\nu(N,h)>-\infty$, then
\begin{equation}\label{eq:positive-avr}
\AVR(N,h)\geq e^{\nu(N,h)}>0.
\end{equation}
\end{lemma}
\begin{proof}
Fix $o\in N$, write $r(x)=d_h(o,x)$, and put
\[
V(s)=\Vol_hB_h(o,s),\qquad
\theta(s)=\frac{V(s)}{\omega_n s^n},\qquad
 a=\AVR(N,h).
\]
For $\tau>0$, define the Gaussian normalization
\begin{equation}\label{eq:I-tau}
I_\tau=(4\pi\tau)^{-n/2}\int_N e^{-r^2/(4\tau)}\dd V_h.
\end{equation}
We first prove that
\begin{equation}\label{eq:I-limit}
\lim_{\tau\to\infty}I_\tau=a.
\end{equation}
This limit also follows from \cite[Lemma~8.10]{ChowEtAl2007}; we give the details here.

The bound $V(s)\leq\omega_n s^n$ implies that $I_\tau$ is finite and strictly positive. More generally, every Gaussian moment
\[
\int_Nr^me^{-cr^2}\dd V_h
\]
is finite for $m\geq0$ and $c>0$. This follows, for example, by decomposing $N$ into the annuli $B_h(o,\ell+1)\setminus B_h(o,\ell)$ and using the Euclidean upper bound for their volumes. These moment bounds will also justify the computations below. The pushforward of $dV_h$ under the distance function is the Stieltjes measure $dV(s)$. Radial integration and integration by parts therefore give
\begin{equation}\label{eq:I-radial}
\begin{split}
I_\tau
 &=(4\pi\tau)^{-n/2}\int_0^\infty e^{-s^2/(4\tau)}\dd V(s)\\
 &=(4\pi\tau)^{-n/2}\int_0^\infty\frac{s}{2\tau}e^{-s^2/(4\tau)}V(s)\dd s\\
 &=\frac{\omega_n}{2(4\pi)^{n/2}}\int_0^\infty
 \theta(\sqrt\tau u)u^{n+1}e^{-u^2/4}\dd u.
\end{split}
\end{equation}
The boundary terms vanish: $V(s)\to0$ as $s\downarrow0$, and the Gaussian dominates the Euclidean upper volume growth as $s\to\infty$. The measure
\[
\frac{\omega_n}{2(4\pi)^{n/2}}u^{n+1}e^{-u^2/4}\dd u
\]
has total mass one. In fact,
\[
\int_0^\infty u^{n+1}e^{-u^2/4}\dd u
=2^{n+1}\Gamma\left(\frac n2+1\right),
\qquad
\omega_n=\frac{\pi^{n/2}}{\Gamma(n/2+1)}.
\]
The substitution $s=\sqrt\tau u$ in \eqref{eq:I-radial} is useful because this probability measure no longer depends on $\tau$. Since $0\leq\theta\leq1$ and $\theta(\sqrt\tau u)\to a$ for every $u>0$, dominated convergence yields \eqref{eq:I-limit}. Moreover, the monotonicity of $\theta$ shows directly from \eqref{eq:I-radial} that $I_\tau$ is nonincreasing as a function of $\tau$.

Consider the normalized Gaussian test function
\begin{equation}\label{eq:gaussian-test}
u_\tau=(4\pi\tau)^{-n/4}I_\tau^{-1/2}e^{-r^2/(8\tau)},
\qquad \int_Nu_\tau^2\dd V_h=1.
\end{equation}
Although $u_\tau$ need not be smooth at the cut locus and does not have compact support, it is admissible by approximation. The distance function is locally Lipschitz, its weak gradient has length one almost everywhere away from $o$, and the cut locus has measure zero. Thus the gradient of $u_\tau$ is well defined almost everywhere. The Gaussian moment bounds give
\[
\int_N\left(u_\tau^2+|\nabla u_\tau|_h^2
+u_\tau^2|\log u_\tau^2|\right)\dd V_h<\infty.
\]
For fixed $\tau$, multiply $u_\tau$ by a radial Lipschitz cutoff $\chi_L$ that equals one on $B_h(o,L)$, vanishes outside $B_h(o,2L)$, and satisfies $0\leq\chi_L\leq1$ and $|d\chi_L|\leq C/L$. Put $w_L=\chi_Lu_\tau$. Then
\[
\norm{w_L-u_\tau}_{L^2}\longrightarrow0,\qquad
\norm{\nabla w_L-\nabla u_\tau}_{L^2}\longrightarrow0.
\]
For the second assertion, use
\[
\nabla w_L-\nabla u_\tau
=(\chi_L-1)\nabla u_\tau+u_\tau\nabla\chi_L.
\]
The first term tends to zero by integrability of $|\nabla u_\tau|^2$, and the squared $L^2$ norm of the second is at most $CL^{-2}\int_Nu_\tau^2\dd V_h$. For the entropy term, write
\[
w_L^2\log w_L^2
=\chi_L^2u_\tau^2\log u_\tau^2+u_\tau^2\chi_L^2\log\chi_L^2.
\]
The first term converges in $L^1$ by dominated convergence. The absolute integral of the second is bounded by $e^{-1}\int_{N\setminus B_h(o,L)}u_\tau^2\dd V_h$, which tends to zero. Finally, $\norm{w_L}_{L^2}\to1$, so replacing $w_L$ by $w_L/\norm{w_L}_{L^2}$ preserves all these limits. Each compactly supported Lipschitz approximation can in turn be approximated uniformly and in $W^{1,2}$ by nonnegative compactly supported smooth functions and normalized again. On a fixed compact set, uniform approximation also gives convergence of the entropy integral, using the continuous extension $t^2\log t^2=0$ at $t=0$. Consequently,
\[
\nu(N,h)\leq\Wcal(h,u_\tau,\tau).
\]

Since $|dr|_h=1$ almost everywhere away from $o$ and $R_h=0$, direct substitution into the $\Wcal$-functional gives
\begin{equation}\label{eq:W-gaussian}
\begin{gathered}
\Wcal(h,u_\tau,\tau)=\frac{J_\tau}{I_\tau}+\log I_\tau-n,\\
J_\tau:=(4\pi\tau)^{-n/2}\int_N\frac{r^2}{2\tau}e^{-r^2/(4\tau)}\dd V_h.
\end{gathered}
\end{equation}
To see the normalization explicitly, $4\tau|\nabla u_\tau|^2=(r^2/(4\tau))u_\tau^2$, whereas
\[
-\log u_\tau^2=\frac n2\log(4\pi\tau)+\log I_\tau+\frac{r^2}{4\tau}.
\]
The two quadratic terms therefore add to the moment in \eqref{eq:W-gaussian}.

On each compact interval of positive $\tau$ values, the differentiated integrand in \eqref{eq:I-tau} is bounded by a fixed polynomial in $r$ times a Gaussian. The moment bounds therefore justify differentiation under the integral sign:
\[
I_\tau'=-\frac n{2\tau}I_\tau
+(4\pi\tau)^{-n/2}\int_N\frac{r^2}{4\tau^2}e^{-r^2/(4\tau)}\dd V_h.
\]
Hence
\begin{equation}\label{eq:I-derivative}
2\tau I_\tau'=J_\tau-nI_\tau,
\qquad \frac{J_\tau}{I_\tau}=n+2\tau(\log I_\tau)'.
\end{equation}
Since $I_\tau$ is nonincreasing, \eqref{eq:W-gaussian} gives
\begin{equation}\label{eq:gaussian-entropy-bound}
\nu(N,h)\leq\Wcal(h,u_\tau,\tau)
=\log I_\tau+2\tau(\log I_\tau)'
\leq\log I_\tau.
\end{equation}
Letting $\tau\to\infty$ and using \eqref{eq:I-limit}, we obtain
\[
\nu(N,h)\leq\log a,
\]
with $\log0=-\infty$. In particular, $a=0$ implies $\nu(N,h)=-\infty$.
\end{proof}

\begin{proposition}[Euclidean product invariance of the volume ratio]\label{prop:product-avr}
Let $(X^d,h_X)$ be a complete connected Ricci-flat manifold, let $k\geq1$, and give $Y=\R^k\times X$ the product metric $h=\gE+h_X$. Then
\begin{equation}\label{eq:product-avr}
\AVR(Y,h)=\AVR(X,h_X).
\end{equation}
\end{proposition}
\begin{proof}
Fix $x_0\in X$ and write
\[
\theta_X(s)=\frac{\Vol_{h_X}B_{h_X}(x_0,s)}{\omega_d s^d}.
\]
The product distance satisfies
\[
d_Y\bigl((z,x),(0,x_0)\bigr)^2=|z|^2+d_X(x,x_0)^2.
\]
Indeed, for a path parametrized on $[0,1]$, the integral of its squared speed in the product is the sum of the corresponding integrals for its two components. Each component integral is bounded below by the square of the corresponding distance. Conversely, the straight segment in $\R^k$ and a minimizing geodesic in $X$, parametrized on the same unit interval with constant speeds, realize this lower bound. Consequently, for fixed $z$ with $|z|<r$, the section of $B_h((0,x_0),r)$ over $z$ is exactly
\[
\{z\}\times B_{h_X}\left(x_0,\sqrt{r^2-|z|^2}\right).
\]
The volume measure of the product is $dz\,dV_{h_X}$. Integrating the corresponding ball sections over the Euclidean factor by Fubini's theorem gives the first equality below. For the second, set $z=ru$, so that $dz=r^kdu$, and use the definition of $\theta_X$. We obtain
\begin{equation}\label{eq:product-volume}
\begin{split}
\Vol_hB_h((0,x_0),r)
 &=\int_{|z|<r}\Vol_{h_X}B_{h_X}\left(x_0,\sqrt{r^2-|z|^2}\right)\dd z\\
 &=\omega_dr^{d+k}\int_{B_1^{\R^k}}(1-|u|^2)^{d/2}
\theta_X\left(r\sqrt{1-|u|^2}\right)\dd u.
\end{split}
\end{equation}
For every $|u|<1$, $\theta_X(r\sqrt{1-|u|^2})\to\AVR(X,h_X)$ as $r\to\infty$. Bishop--Gromov comparison gives $0\leq\theta_X\leq1$, so the integrand is dominated by the integrable function $(1-|u|^2)^{d/2}$. Dividing \eqref{eq:product-volume} by $\omega_{d+k}r^{d+k}$ and applying dominated convergence yields
\[
\AVR(Y,h)=\left(
\frac{\omega_d}{\omega_{d+k}}\int_{B_1^{\R^k}}(1-|u|^2)^{d/2}\dd u
\right)\AVR(X,h_X).
\]
Finally, Fubini's theorem applied to the Euclidean unit ball in $\R^{d+k}$ gives
\begin{equation}\label{eq:euclidean-volume}
\omega_d\int_{B_1^{\R^k}}(1-|u|^2)^{d/2}\dd u=\omega_{d+k}.
\end{equation}
Substitution proves \eqref{eq:product-avr}, including the case in which both asymptotic volume ratios vanish.
\end{proof}

\section{The proof of Theorem~\ref{thm:main}}\label{sec:proof}
We recall the four-dimensional Ricci-flat gap theorem used in the final step.

\begin{proposition}[{\cite[Proposition~3.7]{HondaMondino2026}}]\label{prop:rf-gap}
Let $(X^4,h_X)$ be a complete Ricci-flat four-manifold. If
\[
\AVR(X,h_X)>\frac12,
\]
then $(X,h_X)$ is isometric to Euclidean $\R^4$.
\end{proposition}

\begin{proof}[Proof of Theorem~\ref{thm:main}]
If $M$ is compact, its curvature is bounded by smoothness. Suppose that $M$ is noncompact and $\sup_M|\Rm_g|=\infty$. By Proposition~\ref{prop:blowup} and Lemma~\ref{lem:splitting}, a pointed curvature blow-up produces a complete nonflat Ricci-flat manifold with a global product decomposition
\[
(Y^5,h)=(\R\times X^4,dt^2+h_X).
\]
The factor $X$ is complete and Ricci flat. It is nonflat since the product curvature is that of $X$ and $|\Rm_h|(p_\infty)=1$.

The compactly supported entropy transfer in Proposition~\ref{prop:blowup} and the Gaussian estimate in Lemma~\ref{lem:entropy-avr} give
\begin{equation}\label{eq:final-entropy}
\nu(M,g)\leq\nu(Y,h)\leq\log\AVR(Y,h).
\end{equation}
By \eqref{eq:main-hyp}, the left-hand side is strictly greater than $-\log2$. In particular, $\nu(Y,h)>-\infty$, so Lemma~\ref{lem:entropy-avr} gives positive asymptotic volume ratio. Exponentiating \eqref{eq:final-entropy}, we obtain
\[
\AVR(Y,h)\geq e^{\nu(Y,h)}\geq e^{\nu(M,g)}>\frac12.
\]
Proposition~\ref{prop:product-avr}, with a one-dimensional Euclidean factor, yields
\[
\AVR(X,h_X)=\AVR(Y,h)>\frac12.
\]
Proposition~\ref{prop:rf-gap} now implies $(X,h_X)\cong\R^4$. Consequently $Y$ is flat, contradicting $|\Rm_h|(p_\infty)=1$. Thus the assumption $\sup_M|\Rm_g|=\infty$ is impossible, and the theorem follows.
\end{proof}

\begin{remark}\label{rem:gap-role}
The strict entropy inequality is used only in the final volume-ratio gap. The curvature blow-up and the real-line splitting in Section~\ref{sec:blowup} require only the scalar curvature upper bound for the fixed shrinker. The local estimate in Lemma~\ref{lem:parabolic} provides the Ricci bound needed for the rescaled Hessian calculation.
\end{remark}


\end{document}